\documentclass[11pt]{article}
\usepackage{amscd}
\usepackage{amsfonts}
\usepackage{amsmath}
\usepackage{amssymb}
\usepackage{amsthm}
\usepackage{bbm}
\usepackage{CJK}
\usepackage{fancyhdr}
\usepackage{graphicx}
\usepackage{hyperref}
\usepackage{indentfirst}
\usepackage{latexsym}
\usepackage{mathrsfs}
\usepackage{xypic}

\usepackage[top=1in,bottom=1in,left=1.25in,right=1.25in]{geometry}
\def\<{\langle}
\def\>{\rangle}
\def\a{\alpha}
\def\b{\beta}

\def\D{\Delta}

\def\o{\otimes}

\def\v{\epsilon }

\date{}
\begin{document}
\renewcommand{\baselinestretch}{1.2}
\renewcommand{\arraystretch}{1.0}
\title{\bf The construction and deformation of  BiHom-Novikov agebras}
\author{{\bf Shuangjian Guo$^{1}$, Xiaohui Zhang$^{2}$, Shengxiang Wang$^{3}$\footnote
        { Corresponding author(Shengxiang Wang):~~wangsx-math@163.com} }\\
{\small 1. School of Mathematics and Statistics, Guizhou University of Finance and Economics} \\
{\small  Guiyang  550025, P. R. of China} \\
{\small 2. School of Mathematical Sciences, Qufu Normal University}\\
{\small  Qufu 273165, P. R. of  China}\\
{\small 3.~ School of Mathematics and Finance, Chuzhou University}\\
 {\small   Chuzhou 239000,  P. R. of China}}
 \maketitle
 \maketitle
\begin{center}
\begin{minipage}{13.cm}

{\bf \begin{center} ABSTRACT \end{center}}
BiHom-Novikov agebra is a generalized Hom-Novikov algebra endowed with two commuting multiplicative linear maps. The main purpose of this paper is to   show that two classes of BiHom-Novikov algebras can be
constructed from BiHom-commutative algebras together with derivations and
BiHom-Novikov algebras with Rota-Baxter operators, respectively. We show
that quadratic BiHom-Novikov algebras are associative algebras and
the sub-adjacent BiHom-Lie algebras of BiHom-Novikov algebras are 2-step
nilpotent. Moreover, we develop the 1-parameter formal deformation theory of
BiHom-Novikov algebras


 \end{minipage}
 \end{center}
 \normalsize\vskip1cm

\section*{INTRODUCTION}
\def\theequation{0. \arabic{equation}}
\setcounter{equation} {0}

Novikov algebras were introduced in connection with the Poisson brackets of hydrodynamic
type \cite{Balinskii1985, Dubrovin1983} and Hamiltonian operators in the formal variational calculus \cite{Gelfand1979, Xu1996}. The theoretical study of Novikov algebras was started by Zel¡¯manov and Tikhov \cite{Zelmanov1987} and Filipov \cite{Filipov1989}. But the term ¡°Novikov
algebra¡± was first used by Osborn \cite{Osborn1992}. The left multiplication operators of a Novikov algebra form a
Lie algebra. Thus, it is effective to relate the study of Novikov algebras to the theory of Lie algebras
\cite{Bai2003}. Novikov algebras are a special class of left-symmetric algebras (or under other names such as
pre-Lie algebras, quasi-associative algebras and Vinberg algebras), arising from the study of affine
manifolds, affine structures and convex homogeneous cones \cite{Auslander1977}. Left-symmetric algebras
have close relations with many important fields in mathematics and mathematical physics, such
as infinite-dimensional Lie algebras \cite{Balinskii1985}, classical and quantum Yang-Baxter equation,
quantum field theory and so on.

A quadratic Novikov algebra, introduced in \cite{Bordemann1997}, is a Novikov algebra with a symmetric
nondegenerate invariant bilinear form. The motivation for studying quadratic Novikov algebras
came from the fact that Lie algebras or associative algebras with symmetric nondegenerate invariant
bilinear forms have important applications in several areas of mathematics and physics, such as the
structure theory of finite-dimensional semisimple Lie algebras, the theory of complete integrable
Hamiltonian systems and the classification of statistical models over two-dimensional graphs.

The deformation is a tool to study a mathematical object by deforming it into a family
of the same kind of objects depending on a certain parameter. The deformation theory
was introduced by Gerstenhaber for rings and algebras \cite{Gerstenhaber1964}-\cite{Gerstenhaber1974}, by Kubo and Taniguchi
for Lie triple systems \cite{Kubo2004}, by Ma, Chen and Lin for Hom-Lie Yamaguti algebras \cite{Ma2015} and
Hom-Novikov superalgebras \cite{Sun2015}, by Bai and Meng for Novikov algebras \cite{Bai2001}. They studied 1-
parameter formal deformations and established the connection between the cohomology groups and infinitesimal deformations.

As generalizations of Lie algebras, Hom-Lie algebras were introduced motivated by
applications to physics and to deformations of Lie algebras, especially Lie algebras
of vector fields. The notion of Hom-Lie algebras was firstly introduced by Hartwig,
Larsson and Silvestrov  to describe the structure of certain $q$-deformations
of the Witt and the Virasoro algebras, see \cite{Aizawa, Chaichian, Hartwig, Hu}.
More precisely, a  Hom-Lie algebras are different from
Lie algebras as the Jacobi identity is replaced by a twisted form using a morphism.

The twisting of parts of the defining identities was transferred to other algebraic structures.
In \cite{Makhlouf2008, Makhlouf2009,Makhlouf2010},  Makhlouf and  Silvestrov introduced the notions of
 Hom-associative algebras, Hom-coassociative coalgebras, Hom-bialgebras  and Hom-Hopf algebras.
The original definition of a Hom-bialgebra involved two linear maps, one twisting the
associativity condition and the other one twisting the coassociativity condition.

Yau in \cite{Yau2011} introduced Hom-Novikov algebras, in which the two defining identities are
twisted by a linear map. It turned out that Hom-Novikov algebras can be constructed from
Novikov algebras, commutative Hom-associative algebras and Hom-Lie algebras along
with some suitable linear maps. And Yuan in \cite{Yuan2012} introduced quadratic Hom-
Novikov algebras.  Later, Later, Zhang, Hou and Bai in \cite{Zhang} defined a Hom-Novikov superalgebra
as a twisted generalization of Novikov superalgebras.

A Bihom-algebra is an algebra in such a way that the identities defining the structure
are twisted by two homomorphisms $\a,\b$. This class of algebras was introduced from a
categorical approach in \cite{Graziani} as an extension of the class of Hom-algebras. When the two
linear maps are same automorphisms, Bihom-algebras will be return to Hom-algebras.
These algebraic structures include Bihom-associative algebras, Bihom-Lie algebras and
Bihom-bialgebras. The representation theory of Bihom-Lie algebras was introduced by
Cheng and Qi in \cite{Cheng16}, in which, Bihom-cochain complexes, derivation, central extension, derivation
extension, trivial representation and adjoint representation of Bihom-Lie algebras were
studied. More applications of  BiHom-algebras, BiHom-Lie superalgebras and BiHom-Lie colour algebras
can be found in (\cite{Abdaoui17},\cite{Liu17},\cite{Wang16}).

In the present paper, we consider the construction and deformation of  BiHom-Novikov agebras.

In Section 2, we focus on BiHom-Novikov algebras. We establish the relationships between BiHom-
Novikov algebras and BiHom-Lie algebras as well as the relations between BiHom-Novikov algebras and
Hom-Novikov algebras. Also, we provide a construction of BiHom-Novikov algebras from a commutative
associative algebra  together with derivations and BiHom-Novikov algebras with Rota-Baxter operators, respectively.

In Section 3, we introduce the notion of quadratic BiHom-Novikov algebras and the relationship between BiHom-Novikov
algebras and sub-adjacent BiHom-Lie algebras. We also show quadratic BiHom-
Novikov algebras are associative algebras and the sub-adjacent BiHom-Lie
algebras of BiHom-Novikov algebras are 2-step nilpotent.

In Section 4, we define low orders coboundary operator and give low order cohomology groups
of BiHom-Novikov algebras. We show that the cohomology group is suitable for this
1-parameter formal deformation theory.

\section{Preliminaries}
\def\theequation{\arabic{section}.\arabic{equation}}
\setcounter{equation} {0}

Throughout the paper, all algebraic systems are supposed to be over a field $\mathbb{K}$.
 Any unexplained definitions and notations can be found in \cite{Graziani} and \cite{Yau2010}.

In this section we recall some basic definitions and results related to our paper.

\noindent{\bf 1.1. BiHom-associative algebra} A  BiHom-associative algebra is a 4-tuple $(A,\mu,\alpha,\beta)$,
where $A$ is a ${k}$-linear space, $\alpha: A\rightarrow A$, $\beta: A\rightarrow A$ and
 $\mu: A\o A\rightarrow A$  are linear maps, with notation $\mu(a\o b)=ab$,
 satisfying the following conditions, for all  $a,a',a''\in A$:
\begin{eqnarray*}
&&\alpha\circ\beta=\beta\circ\alpha,\\
&&\alpha(aa')=\alpha(a)\alpha(a'),\beta(aa')=\beta(a)\beta(a'),~~~~~~~~~~~~~~~~~~~~~~~~~~~~~~~~~~~~~~~~~~~~~~~~~~~~~~~\\
&&\alpha(a)(a'a'')=(aa')\beta(a'').
\end{eqnarray*}
And the maps $\alpha, \beta$ are called the structure maps of $A$.
\smallskip

Clearly, a Hom-associative algebra $(A,\mu,\alpha)$ can be regarded as the BiHom-associative
algebra $(A,\mu,\alpha,\alpha)$.
\medskip

\noindent{\bf 1.2. BiHom-Lie algebra }  A BiHom-Lie algebra
is a 4-tuple $(L,[\cdot,\cdot],\alpha,\beta)$,
where $L$ is a ${k}$-linear space, $\alpha: L\rightarrow L$, $\beta: L\rightarrow L$ and
 $[\cdot,\cdot]: L\o L\rightarrow L$ are linear maps,
 satisfying the following conditions, for all  $a,a',a''\in A$:
\begin{eqnarray*}
&&\alpha\circ\beta=\beta\circ\alpha,\\
&&\alpha[a,a']=[\alpha(a),\alpha(a')],\beta[a,a']=[\beta(a),\beta(a')],\\
&&[\beta(a),\alpha(a')]=-[\beta(a'),\alpha(a)].\\
&&[\beta^{2}(a),[\beta(a'),\alpha(a'')]]+[\beta^{2}(a'),[\beta(a''),\alpha(a)]]+[\beta^{2}(a''),[\beta(a),\alpha(a')]]=0.~~~~~~~~~~~~~~
\end{eqnarray*}

Obviously, a Hom-Lie algebra $(L,[\cdot,\cdot],\alpha)$ is a particular case of a BiHom-Lie
algebra, namely $(L,[\cdot,\cdot],\alpha,\alpha)$.
 Conversely, a BiHom-Lie algebra $(L,[\cdot,\cdot],\alpha,\alpha)$ with bijective $\alpha$
is the Hom-Lie algebra $(L,[\cdot,\cdot],\alpha)$ .

\noindent{\bf 1.3. Hom-Novikov algebra }
A Hom-Novikov algebra is a triple $(A,\mu,\alpha)$ consisting of a vector space $A$, $\alpha: A\rightarrow A$ is a algebra homomorphisms, $\mu: A\o A\rightarrow A$  is an even bilinear map, with notation $\mu(a\o b)=ab$ satisfying
\begin{eqnarray*}
&& (xy)\a(z)=(xz)\a(y),\\
&& (xy)\a(z)-\a(x)(yz)=(yx)\a(z)-\a(y)(xz),
\end{eqnarray*}
for all $x,y,z\in A$.

\section{BiHom-Novikov algebras}
\def\theequation{\arabic{section}. \arabic{equation}}
\setcounter{equation} {0}

 In this section, we establish the relationships between BiHom-Novikov algebras and Hom-Novikov
algebras as well as the relations between BiHom-Novikov algebras and BiHom-Lie algebras. Also, we
provide a construction of BiHom-Novikov algebras from a commutative BiHom-associative algebra along
with a suitable linear map but not a derivation.

\noindent{\bf Definition 2.1.}
A BiHom-Novikov algebra is a 4-tuple $(A,\mu,\alpha,\beta)$ consisting of a vector space $A$, $\alpha: A\rightarrow A$ and $\beta: A\rightarrow A$ are algebras homomorphisms, $\mu: A\o A\rightarrow A$  is an even bilinear map, with notation $\mu(a\o b)=ab$ satisfying
\begin{eqnarray}
&&\alpha\circ\beta=\beta\circ\alpha,\\
&& (xy)\a(z)=(xz)\a(y),\\
&& (\b(x)\a(y))\b(z)-\a\b(x)(\a(y)z)=(\b(y)\a(x))\b(z)-\a\b(y)(\a(x)z),
\end{eqnarray}
for all$x,y,z\in A$.

Clearly, Hom-Novikov algebras and Novikov algebras   are examples of BiHom-Novikov algebras by setting $\b=\a$ and $\a=\b=id$, where
$id$ is the identity map. It was shown in \cite{Yau2010} that any Novikov algebra can be deformed into a
Hom-Novikov algebra along with an algebra endomorphism.

\noindent{\bf Proposition  2.2.}
Let $(A, \mu, \a,\b)$ be a BiHom-Novikov algebra, for all $x,y,z\in A$, we have
\begin{eqnarray}
[\b(x),\a(y)]\a^2(z)+[\b(y),\a(z)]\a^{2}(x)+[\b(z),\a(x)]\a^{2}(y)=0,\\
\b^2(z)[\b(x),\a(y)]+\b^{2}(x)[\b(y),\a(z)]+\b^{2}(y)[\b(z),\a(x)]=0,
\end{eqnarray}
where $\mu(x,y)=xy$ and $[\b(x),\a(y)]=\b(x)\a(y)-\b(x)\a(y)$.

{\bf Proof.} For any $x,y,z\in A$, we have
\begin{eqnarray*}
&&[\b(x),\a(y)]\a^2(z)=(\b(x)\a(y))\a^2(z)-(\b(y)\a(x))\a^2(z),\\
&&[\b(y),\a(z)]\a^2(x)=(\b(y)\a(z))\a^2(x)-(\b(z)\a(y))\a^2(x),\\
&&[\b(z),\a(x)]\a^2(y)=(\b(z)\a(x))\a^2(y)-(\b(x)\a(z))\a^2(y),
\end{eqnarray*}
Since equation (2.2) holds for all$x,y,z\in A$, we have
\begin{eqnarray*}
&&[\b(x),\a(y)]\a^2(z)+[\b(y),\a(z)]\a^{2}(x)+[\b(z),\a(x)]\a^{2}(y)\\
&=&(\b(x)\a(y))\a^2(z)-(\b(y)\a(x))\a^2(z)+(\b(y)\a(z))\a^2(x)-(\b(z)\a(y))\a^2(x)\\
&&+(\b(z)\a(x))\a^2(y)-(\b(x)\a(z))\a^2(y)\\
&=&0,
\end{eqnarray*}
which proves equation (2.4). Now using equations (2.2) and (2.4), we have
\begin{eqnarray*}
&&\b^2(z)[\b(x),\a(y)]+\b^{2}(x)[\b(y),\a(z)]+\b^{2}(y)[\b(z),\a(x)]\\
&=& \b^2(z)(\b(x)\a(y))+\b^{2}(x)(\b(y)\a(z))-\b^{2}(x)(\b(z)\a(y))\\
&&+\b^{2}(y)(\b(z)\a(x))-\b^{2}(y)(\b(x)\a(z))\\
&=&\b^2(z)(\b(x)\a(y))-\b^{2}(x)(\b(z)\a(y))+\b^{2}(x)(\b(y)\a(z))-\b^{2}(y)(\b(x)\a(z))\\
&&\b^{2}(y)(\b(z)\a(x))-\b^2(z)(\b(y)\a(x))\\
&=&(\a^{-1}\b^2(z)\b(x)-\a^{-1}\b^2(x)\b(z))\a\b(y)+(\a^{-1}\b^2(x)\b(y)-\a^{-1}\b^2(y)\b(x))\a\b(z)\\
&&(\a^{-1}\b^2(x)\b(y)-\a^{-1}\b^2(y)\b(x))\a\b(z)\\
&=&[\a^{-1}\b^2(z), \b(x)]\a\b(y)+[\a^{-1}\b^2(y), \b(z)]\a\b(x)+[\a^{-1}\b^2(x), \b(y)]\a\b(z)\\
&=&\a^{-1}\b([\b(x),\a(y)]\a^2(z)+[\b(y),\a(z)]\a^{2}(x)+[\b(z),\a(x)]\a^{2}(y))\\
&=&0.
\end{eqnarray*}
And this finishes the proof.  \hfill $\square$

\noindent{\bf Corollary 2.3.} Let $(A, \mu, \a,\b)$ be a BiHom-Novikov algebra, Define a bilinear map $[\cdot, \cdot]: A\times A\rightarrow A$ by
\begin{eqnarray*}
[x,y]=\mu(x\o y)-\mu(\a^{-1}\b(y)\o \a\b^{-1}x),  \mbox{for any $x,y\in A$}.
\end{eqnarray*}
Then $(A,[\cdot, \cdot],\a,\b )$ is a BiHom-Lie algerba.

\noindent{\bf Definition 2.4.}  Let $(A, \mu, \a,\b)$ be a BiHom-Novikov algebra, Which is called

(i) regular  if $\a,\b$ are algebra automorphism.

(ii) involutive if $\a,\b$ are involution, i.e., $\a^2=\b^2=id$.

Yau \cite{Yau2011} and Yuan \cite{Yuan2012} gave a way to construct Hom-Novikov algebras, starting from a Novikov algebra
and an algebra endomorphism. In the following, we provide a construction of BiHom-Novikov algebras
from Hom-Novikov algebras and Novikov algebras  along with algebra automorphisms.

\noindent{\bf Proposition  2.5.} Let $(A, \mu, \a,\b)$ be  an involutive   BiHom-Novikov algebra,   define a new multiplication on $A$ by $x\ast y=\a(x)\b(y)$.  Then $(A, \ast)$ is a Novikov algebra.

{\bf Proof.}  We write $\mu(x\o y)=xy$ and $x\ast y=\b(x)\a(y)$, for all $x,y\in A$. Hence, it needs to show
\begin{eqnarray}
(x\ast y)\ast z=(x\ast z)\ast y\\
(x\ast y) \ast z-x\ast(y \ast z)=(y \ast x)\ast z-y\ast (x\ast z),
\end{eqnarray}
for all $x,y,z\in A$. Since $\a,\b$ are involution, we have
\begin{eqnarray*}
(x\ast y)\ast z=\a(\a(x)\b(y))\b(z)=(\a^2(x)\a\b(y))\b(z)=
\end{eqnarray*}
Similar, we have $(x\ast z)\ast y=\a(\a(x)\b(z))\b(y)=(\a^2(x)\a\b(z))\b(y)$, from which equation (2.6) follows since  $(A, \mu, \a,\b)$ is  an involutive   BiHom-Novikov algebra. Furthermore, using equation (2.3), we have
\begin{eqnarray*}
&&(x\ast y)\ast z-x\ast(y\ast z)\\
&=&\a(\a(x)\v(y))\b(z)-\a(x)\b(\a(y)\b(z))\\
&=& (\a^2(x)\a\b(y))\b(z)-\a(x)(\a\b(y)\b^2(z))\\
&=&(\b^2(y)\a\b^{-1}(x))\b(z)-\a(y)(\a\b(x)z)\\
&=&\a(\a(y)\b(x))\b(z)-\a(y)\b(\a(x)\b(z))\\
&=&(y\ast x)\ast z-y\ast(x\ast z),
\end{eqnarray*}
which proves equation (2.7) and the proposition.   \hfill $\square$

\noindent{\bf Proposition  2.6.}  Let $(A, \mu)$ be  a  Novikov algebra, and $\a,\b: A\rightarrow A $ are linear maps satisfying $\a\circ \b=\b\circ \a, \a(x\o y)=\a(x)\a(y)$ and $\b(xy))=\b(x)\b(y)$, for all $x,y \in A$, and we  define a new multiplication on $A$ by $x\star y=\a(x)\b(y)$.  Then $(A, \star, \a,\b)$ is a BiHom-Novikov algebra.

{\bf Proof.}  We only prove (2.3) and leave the rest to the reader: for any $x,y,z\in A$, we have
\begin{eqnarray*}
 && (\b(x)\star\a(y))\star\b(z)-\a\b(x)\star(\a(y)\star z)\\
 &=& \a(\a\b(x)\a\b(y))\b^2(z)-\a^2\b(x)(\a^2\b(y)\b^2(z))\\
 &=&(\a^2\b(y)\a^2\b(x))\b^2(z)-\a^2\b(y)(\a^2\b(x)\b^2(z))\\
&=& \a(\a\b(y)\a\b(x))\b^2(z)-\a^2\b(y)(\a^2\b(x)\b^2(z))\\
 &=&(\b(y)\star\a(x))\star\b(z)-\a\b(y)\star(\a(x)\star z)
\end{eqnarray*}

\noindent{\bf Proposition  2.7.}
Let $(A, \mu, \a,\b)$ be  a regular   BiHom-Novikov algebra,   define a new multiplication on $A$ by $[x,y]'=[\a^{-1}(x),\b^{-1}(y)]$
where $[x,y]=xy-\alpha^{-1}(\beta(y))\alpha(\beta^{-1}(x))$, for all $x,y\in A$, then $(A, [\cdot,\cdot]')$ is a Lie algebra. In particular, if $\a,\b$ are involution and we define multiplication on $A$ by $[x,y]''=[\a(x),\b(y)]$, then $(A, [\cdot,\cdot]'')$ is a Lie algebra.

{\bf Proof.} For any $x,y,z\in A$, we have
\begin{eqnarray*}
&&[[x,y]',z]'\\
&=&[[\a^{-2}(x),\a^{-1}\b^{-1}(y)],\b^{-1}(z)]\\
&=&[\a^{-2}(x),\a^{-1}\b^{-1}(y)]\b^{-1}(z)-\a^{-1}(z)\a\b^{-1}[\a^{-2}(x),\a^{-1}\b^{-1}(y)]\\
&=&\a^{-2}\b^{-1}([\b(x),\a(y)]\a^{2}(z)-\a\b(z)[\a(x),\a^{2}\b^{-1}(y)])
\end{eqnarray*}
Similarly, we have
\begin{eqnarray*}
&&[[y,z]',x]'=\a^{-2}\b^{-1}([\b(y),\a(z)]\a^{2}(x)-\a\b(x)[\a(y),\a^{2}\b^{-1}(z)]),\\
&&[[z,x]',y]'=\a^{-2}\b^{-1}([\b(z),\a(x)]\a^{2}(y)-\a\b(y)[\a(z),\a^{2}\b^{-1}(x)]),
\end{eqnarray*}
Then it follows from equations (2.4) and (2.5) that
\begin{eqnarray*}
[[x,y]',z]'+[[y,z]',x]'+[[z,x]',y]'=0
\end{eqnarray*}
Clearly, $[\b(x),\a(y)]=-[\b(y),\a(x)]$, So we have
\begin{eqnarray*}
[x,y]'=-[y,x]'
\end{eqnarray*}
It follows immediately that $(A, [\cdot,\cdot]'')$ is also a Lie algebra  when $\a,\b$ are involution.  \hfill $\square$

Let $(A, \mu)$ be a commutative associative algebra,  $\a,\b: A\rightarrow A$ be two commuting algebra morphism,   and  $D: A\rightarrow A$ be an even derivation with $D\a=\a D, D\b=\b A$. Consider the following operation on $A$:
\begin{eqnarray}
x\star y=\a(x)D(\b(y)), ~~~\mbox{for any $x,y\in A$}.
\end{eqnarray}
Next we check equation (2.2), on one hand,
\begin{eqnarray*}
(x\star y)\star \a(z)&=& (\a(x)\D(\b(y)))\star \a(z)\\
&=&  (\a^{2}(x)D(\a\b(y)))D(\a\b(z))\\
&=& \a^2(x)(D(\a\b(y))D(\a\b(z))),
\end{eqnarray*}
on the other hand,
\begin{eqnarray*}
(x\star z)\star\a(y)&=& (\a(x)\D(\b(z)))\star \a(y)\\
&=&  (\a^{2}(x)D(\a\b(z)))D(\a\b(y))\\
&=& \a^2(x)(D(\a\b(z))D(\a\b(y)))\\
&=&\a^2(x)(D(\a\b(y))D(\a\b(z))),
\end{eqnarray*}
hence£¬ we have
\begin{eqnarray*}
(x\star y)\star \a(z)=(x\star z)\star\a(y).
\end{eqnarray*}
Similarly, we have
\begin{eqnarray*}
&&(\b(x)\star\a(y))\star\b(z)-\a\b(x)(\a(y)z)\\
&=&\a^2\b(x)(D(\a^2\b(y))D(\b^2(z)))-\a^2\b(x)(D(\a^2\b(y)D(\b(z))),
\end{eqnarray*}
and
\begin{eqnarray*}
&&(\b(y)\star\a(x))\star\b(z)-\a\b(y)(\a(x)z)\\
&=&\a^2\b(y)(D\a^2\b(x)D(\b^2(z)))-\a^2\b(y)(D(\a^2\b(x)D(\b(z))).
\end{eqnarray*}
\hfill $\square$

Now from the discussions above, we obtain

\noindent{\bf Theorem 2.8.} Let $(A, \mu)$ be a commutative associative algebra,  $\a,\b: A\rightarrow A$ be two commuting algebra morphism,   and  $D: A\rightarrow A$ be an even derivation with $D\a=\a D, D\b=\b A$. Then $(A, \mu, \a,\b)$ is a BiHom-Novikov algebra, where
$\star$ is defined by (2.8).

\noindent{\bf Definition  2.9.} Let $(A, \star, \a,\b)$ be a BiHom-algebra and let $\lambda\in k$. If a linear map $P: A\rightarrow A$ stisfies
\begin{eqnarray*}
P(x)\star P(y)=P(P(x)\star y+x\star P(y)+\lambda x\star y), ~~~\mbox{for any $x,y\in A$}.
\end{eqnarray*}
Then $P$ is called  a Rota-Baxter operator of weight $\lambda$ and $(A, \star, \a,\b, P)$ is called a Rota-Baxter
Hom-algebra of weight $\lambda$.

\noindent{\bf Theorem 2.10.}  Let $(A, \star, \a,\b, P)$ be a Rota-Baxter BiHom-Novikov algebra of weight
$\lambda$ and $P$ an even linear map. Assume that $\a, \b$ and $P$ commute,  define a new multiplication on $A$ by
\begin{eqnarray*}
x\circ y=P(x)\star y+x\star P(y)+\lambda x\star y,~~~\mbox{for any $x,y\in A$}.
\end{eqnarray*}
 Then $(A, \circ, \a,\b, P)$ is a  BiHom-Novikov algebra.

{\bf Proof.} The multiplication of $\a,\b$ with respect to $\circ$ follows from the multiplication of
$\a,\b$ with respect to $\star$ and the hypothesis $P\a=\a P, P\b=\b P$. We only check  equation (2.2) and (2.3 is left to the reader.) For any $x,y,z\in A$, on one hand,  we have,
\begin{eqnarray*}
&&(x\circ y) \circ\a(z)\\
&=& P(P(x)\star y+x\star P(y)+\lambda x\star y)\star \a(z)+(P(x)\star y+x\star P(y)+\lambda x\star y)\star P(\a(z))\\
&&+\lambda(P(x)\star y+x\star P(y)+\lambda x\star y)\star \a(z)\\
&=& (P(x)\star P(y))\star\a(z)+(P(x)\star y)\star \a(P(z))+(x\star P(y))\star\a(P(z))+\lambda(x\star y)\star\a(P(z))\\
&&\lambda(P(x)\star y)\star\a(z)+\lambda(x\star P(y))\star\a(z)+\lambda^2(x\star y)\star\a(z)\\
&=&(P(x)\star z)\star\a(P(y))+(P(x)\star P(z))\star \a(P(y))+(x\star P(z))\star\a(P(y))+\lambda(x\star P(z))\star\a(y)\\
&&\lambda(P(x)\star z)\star\a(y)+\lambda(x\star z)\star\a(P(y))+\lambda^2(x\star z)\star\a(y).
\end{eqnarray*}
On the other hand, we have
\begin{eqnarray*}
&&(x\circ z)\circ \a(y)\\
&=&(P(x)\star z)\star\a(P(y))+(P(x)\star P(z))\star \a(P(y))+(x\star P(z))\star\a(P(y))+\lambda(x\star P(z))\star\a(y)\\
&&\lambda(P(x)\star z)\star\a(y)+\lambda(x\star z)\star\a(P(y))+\lambda^2(x\star z)\star\a(y).
\end{eqnarray*}
Hence, the conclusion holds.
\hfill $\square$
\section{Quadratic BiHom-Novikov algebras}
\def\theequation{\arabic{section}. \arabic{equation}}
\setcounter{equation} {0}

In this section, we extend the notions of quadratic Novikov algebras   and quadratic Hom-Novikov algebras to quadratic BiHom-Novikov
algebras and provide some properties. we always assume that the structure maps $\alpha$ and $\beta$ are bijective.

\noindent{\bf Definition 3.1.} Let $(A,[\cdot, \cdot],\a,\b)$ be a BiHom-Lie algebra and $B: A\times A\rightarrow k$ be a  bilinear form  on $A$.

(1) $B$ is said  nondegenerate if
\begin{eqnarray*}
A^{\perp}=\{x\in A|B(x,y)=0,\forall y\in A\}=0.
\end{eqnarray*}

(2)$B$ is said  symmetric if
\begin{eqnarray*}
B(x,y)=B(y,x),~~\mbox{for any $x,y\in A$}.
\end{eqnarray*}

(3)$B$ is said $\a\b$-invariant if
\begin{eqnarray*}
B([\b(x),\a(y)],\a(z))=B(\a(x),[\b(y), \a(z)]), ~~\mbox{for any $x,y,z\in A$}.
\end{eqnarray*}

\noindent{\bf Definition 3.2.}
A quadratic BiHom-Lie algebra is a quintuple $(A, [\cdot,\cdot], \a, \b, B)$ such that $(A, [\cdot,\cdot], \a, \b)$
is a BiHom-Lie algebra with a symmetric invariant nondegenerate bilinear form $B$ satisfying
\begin{eqnarray}
B(\a(x), y)=B(x, \a(y)), B(\b(x), y)=B(x, \b(y))~~~\mbox{for any $x,y\in A$}.
\end{eqnarray}

We can define quadratic BiHom-Novikov algebras as follows.

\noindent{\bf Definition 3.2.} A quadratic BiHom-Novikov algebra is a quintuple $(A, \mu, \a, \b, B)$ is a BiHom-Novikov algebra $(A, \mu, \a, \b)$
 with a symmetric invariant nondegenerate bilinear form $B$ satisfying
\begin{eqnarray}
B(\a(x), \b(y)\a(z))=B(\b(x)\a(y), \a(z)), ~~~\mbox{for any $x,y, z\in A$}.
\end{eqnarray}

 \noindent{\bf Proposition 3.3.}
Let $(A, \mu, \a, \b, B)$ be a quadratic BiHom-Novikov algebra and $(A, [\cdot,\cdot], \a, \b)$
be the sub-adjacent  BiHom-Lie algebra of $A$, if
\begin{eqnarray}
B(\a(x), y)=B(x, \a(y)), B(\b(x), y)=B(x, \b(y)),~~~\mbox{for any $x,y\in A$}.
\end{eqnarray}
Then $(A, [\cdot,\cdot], \a, \b, B_{\a,\b})$  is a quadratic BiHom-Lie algebra, where $B_{\a,\b}(x,y)=B(\a(x),y)$.

{\bf Proof.} Since $B$ is a nondegenerate bilinear form, $B_\a$ is a nondegenerate
bilinear form on $A$. For all $x,y,z\in A$, using the properties of $B$, we have
\begin{eqnarray*}
B_{\a}([\b(x),\a(y)],\a(z))&=&B(\a[\b(x),\a(y)], \a(z))\\
&=&B([\b(x),\a(y)], \a^2(z))\\
&=&B(\b(x)\a(y), \a^2(z))-B(\b(y)\a(x),\a^2(z))\\
&=&B(\a^2(x),\b(y)\a(z))-B(\a^2(x), \b(z)\a(y))\\
&=&B(\a^2(x),[\b(y),\a(z)])\\
&=& B_{\a}(\a(x),[\b(y),\a(z)]).
\end{eqnarray*}
Hence, $B_{\a}$ is invariant. Using symmetry of $B$ and equation (3.3), we have
\begin{eqnarray*}
B_{\a}(x,y)=B(\a(x),y)=B(y,\a(x))=B(\a(y),x)=B_{\a}(y,x).
\end{eqnarray*}
which proves $B_{\a}$ is symmetric. Using equation (3.3) again, we have
\begin{eqnarray*}
B_\a(\a(x),y)=B(\a(\a(x)), y)=B(\a(x),\a(y))=B_\a(x,\a(y)),\\
B_\a(\b(x),y)=B(\a(\b(x)),y)=B(\a(x),\b(y))=B_\a(x,\b(y)).
\end{eqnarray*}
\hfill $\square$

 \noindent{\bf Corollary 3.4.} Let $(A, \mu, \a, \b, B)$ be a quadratic BiHom-Novikov algebra and $\a,\b$ satisfying
equation (3.4) and $(A, [\cdot,\cdot])$
be the sub-adjacent Lie algebra,  define a new multiplication on $A$ by $[x,y]'=[\a(x),\b(y)]$. Then $(A, [\cdot,\cdot]',\a,\b,B_\a)$
is  a quadratic Hom-Lie algebra, where $B_{\a,\b}(x,y)=B(\a(x),y)$.

{\bf Proof.}  Using the similar
arguments as those in the proof of Proposition 3.3, we get $B_{\a}$ is a symmetric nondegenerate bilinear
form with equation (3.1) satisfied. It remains to show that $B_{\a}$ is invariant. For all $x,y,z\in A$,
using invariance and symmetry of $B$, we have
\begin{eqnarray*}
B_{\a}([\b(x),\a(y)]',\a(z))&=&B(\a[\a\b(x),\a\b(y)], \a(z))\\
&=&B([\a\b(x),\a\b(y)], \a^2(z))\\
&=&B(\a\b(x)\a\b(y), \a^2(z))-B(\b^2(y)\a^2(x),\a^2(z))\\
&=&B(\a^2(x),\a\b(y)\a\b(z))-B(\a^2(x), \b(z)\a(y))\\
&=&B(\a^2(x),[\a\b(y),\a\b(z)])\\
&=& B_{\a}(\a(x),[\b(y),\a(z)]').
\end{eqnarray*}
\hfill $\square$

 \noindent{\bf Proposition 3.5.} Let $(A, \mu, \a,\b, B)$ be a quadratic BiHom-Novikov algebra, where $\a,\b$ are involution
satisfying equation (3.3), define a new multiplication on $A$ by $x\star y=\a(x)\star\b(y)$  Then  $(A, \star, \a,\b, B)$ is a quadratic Novikov algebra.

{\bf Proof.} $(A, \mu, \a,\b, B)$ is a BiHom-Novikov algebra by Proposition 2.5. It suffices to show that $B$
is invariant under the operation $\star$. For all $x,y,z\in A$, we have
\begin{eqnarray*}
B(\a(x), \b(y)\star\a(z))=B(\a(x), \a\b(y)\a\b(z))\\
=B(\a\b(x)\a\b(y), \a(z))=B(\b(x)\star\a(y), \a(z)),
\end{eqnarray*}
which completes the proof.
\hfill $\square$

Let $(A, \mu, \a,\b)$ be a  BiHom-Novikov algebra, whose center is denoted by $Z(A)$ and
defined by
\begin{eqnarray*}
Z(A)=\{x\in A|xy=yx=0,\forall y\in A\}.
\end{eqnarray*}
Let  $(A, [\cdot, \cdot], \a',\b')$ be a BiHom-Lie algebra. The lower central series of $A$ is defined as
usual, i.e., $A^{0}= A, A^{i} = [A,A^{i-1}], \forall i \geq 1$. We call $A$ is $i$-step nilpotent if $A^{i} = 0$ and
$A^{i-1}\neq 0$. The center of the Hom-Lie superalgebra is denoted by $C(A)$ and defined by
\begin{eqnarray*}
C(A)=\{x\in A|[x,y]=0,\forall y\in A\}
\end{eqnarray*}

\noindent{\bf Theorem  3.6.}   Let $(A, \mu, \a,\b, B)$ be a quadratic BiHom-Novikov algebra, and $HLie(A)$ be the sub-adjacent BiHom-Lie algebra. Then
$[\b(x),\a(y)]\in Z(A)$, for any $x,y \in HLie(A)$. As a consequence, $HLie(A)$ is 2-step nilpotent.

{\bf Proof.} For any $x,y,z\in A$, we have
\begin{eqnarray*}
\a(x)[y',z']=\a(x)[\b(y),\a(z)]&=&\a(x)(\b(y)\a(z))-\a(x)(\b(y)\a(z))\\
&=&(x\b(y))\a\b(z)-(x\b(y)\a\b(z)\\
&=&(x\b(y))\a\b(z)-(x\b(y)\a\b(z)\\
&=&0.
\end{eqnarray*}
Using Equation (3.2), we have
\begin{eqnarray*}
B([y',z']\a(d), \a^2(x))=B(\a(x)\a\b^{-1}[y',z'], \a^2(d))=0
\end{eqnarray*}
which implies $[x, y]\in Z(A)$ since $\a,\b$ are automorphism and $B$ is nondegenerate. Hence,
we have $[HLie(A),HLie(A)] \in Z(A)$. Obviously, $Z(A) \in C(HLie(A))$. Then it follows
that $HLie(A)$ is 2-step nilpotent.

\section{1-parameter formal deformations of BiHom-Novikov algebras}
\def\theequation{\arabic{section}. \arabic{equation}}
\setcounter{equation} {0}

 \noindent{\bf Definition 4.1.}
Let $(A,\ast, \a,\b)$ be a regular BiHom-Novikov algebra. If an $n$-linear map $f: A\times \cdot \cdot \cdot \times A\rightarrow A$ satisfies
\begin{eqnarray}
\a(f(x_1,\cdot\cdot\cdot,x_n))=f(\a(x_1),\cdot\cdot\cdot,\a(x_n)),\nonumber \\
\b(f(x_1,\cdot\cdot\cdot,x_n))=f(\b(x_1),\cdot\cdot\cdot,\b(x_n)),
\end{eqnarray}
then $f$ is called an $n$-BiHom-cochain on $A$. Denote by $C_{\a,\b}^n(A,A)$ the set of all $n$-Hom-
cochains, $\forall n\geq 1$.

 \noindent{\bf Definition 4.2.} For $n=1,2$, the coboundary operator $\delta_{hom}^n:C_{\a,\b}^n(A,A)\rightarrow C_{\a,\b}^{n+1}(A,A)$ is defined as follows
 \begin{eqnarray*}
 &&\delta^{1}_{hom} f(x_1,x_2)=x_1\ast f(x_2)+f(x_1)\ast x_2 -f(\a^{-1}\b(x_1)\ast x_2),\\
 &&\delta^{2}_{hom} f(x_1,x_2, x_3)=f(\b(x_1), \a^{-1}\b(x_2) \ast x_3)- f(\a^{-1}\b(x_1)\ast x_2, \b(x_3)) \nonumber\\
 &&\hspace{3cm} -f(\b(x_2),\a^{-1}\b(x_1) \ast x_3)+f(\a^{-1}\b(x_2)\ast x_1, \b(x_3)) \\
 &&\hspace{3cm}+\b(x_1)\ast f(x_2,x_3-f(x_1, x_2)\ast \b(x_3) \\
 &&\hspace{3cm}+\b(x_2)\ast f(x_1, x_3)+f(x_2,x_1)\ast \b(x_3).
 \end{eqnarray*}

It is not difficult to verify that $\delta_{hom}^{1},  \delta_{hom}^{2}$  satisfies (4.1). Thus, the coboundary
operator $\delta_{hom}^{n}$  is well-defined.

\noindent{\bf Theorem 4.3.} The coboundary operator$\delta_{hom}^{1},  \delta_{hom}^{2}$  defined above satisfies $\delta_{hom}^{2}\delta_{hom}^{1}=0$.

{\bf Proof.} Similar to \cite{Sun2015} and left to the reader.
\hfill $\square$

For $n = 1, 2$, the map $f\in C_{\a,\b}^n(A,A)$ is called an $n$-Hom-cocycle if $\delta^{n}_{hom} f=0$. We denote $Z_{\a,\b}^n(A,A)$ the subspace spanned by $n$-BiHom-cocycles and $B_{\a,\b}^n(A,A)=\delta^{n-1}_{hom}C_{\a,\b}^{n-1}(A,A)$. Since $\delta_{hom}^{2}\delta_{hom}^{1}=0$, $B_{\a,\b}^n(A,A)$
 is a subspace of $Z_{\a,\b}^n(A,A)$. Hence we can define a cohomology space $H_{\a,\b}^2(A,A)$  of $(A,\ast, \a,\b)$  as the factor space $Z_{\a,\b}^2(A,A)/B_{\a,\b}^2(A,A)$.

Let $(A,\ast, \a,\b)$ be a regular BiHom-Novikov algebra and $K[[t]]$ be the ring of formal power
series over $\mathbb{K}$. Suppose that $A[[t]]$ is the set of formal power series over $A$. Then for an
$\mathbb{K}$-bilinear map $f:A\times A\rightarrow A$, it is natural to extend it to be a $\mathbb{K}[[t]]$-bilinear map
$f:A[[t]]\times A[[t]]\rightarrow A[[t]]$ by
\begin{eqnarray*}
f(\sum_{i\geq 0}x_it^i, \sum_{j\geq 0}y_jt^j)=\sum_{i,j\geq 0}f(x_i,y_j)t^{i+j}.
\end{eqnarray*}

 \noindent{\bf Definition 4.4.} Let $(A,\ast, \a,\b)$ be a regular  BiHom-Novikov algebra over $\mathbb{K}$. A 1-parameter
formal deformation of $(A,\ast, \a,\b)$ is a formal power series $g_t:A[[t]]\times A[[t]]\rightarrow A[[t]]$ of the
form
\begin{eqnarray*}
g_t(x,y)=\sum_{i\geq 0} G_i(x,y)t^i= G_0(x,y)+ G_1(x,y)t+ G_2(x,y)t^2+\cdot\cdot\cdot,
\end{eqnarray*}
where each $G_i$ is an even $\mathbb{K}$-bilinear map $G_i: A\times A\rightarrow A$ and $G_0(x,y)=x\ast y$ such that the following identities hold
\begin{eqnarray}
&&g_t(\a(x),\a(y))=\a\circ g_t(x,y), \\
&&g_t(\b(x),\b(y))=\b\circ g_t(x,y),\\
&&g_t(g_t(x,y), \a(z))=g_t(g_t(x,z),\a(y)),\\
&&g_t(\a\b(x),g_t(\a(y),z))-g_t(g_t(\b(x)\a(y)),\b(z))\nonumber\\
&&- g_t(\a\b(y),g_t(\a(x), z))+g_t(g_t(\b(y),\a(x)), \b(z))=0.
\end{eqnarray}
Conditions (4.2)-(4.5) are called the deformation equations of a regular BiHom-Novikov algebra.

Note that $A[[t]]$ is a module over $K[[t]]$ and $g_t$ defines the bilinear multiplication on
$A[[t]]$ such that $A_t = (A[[t]], g_t, \a,\b)$ is a BiHom-Novikov algebra. Now we investigate
the deformation equations (4.2)-(4.5).

Conditions (4.2)-(4.5) are equivalent to the following equations
\begin{eqnarray}
&&G_i(\a(x),\a(y))=\a\circ G_i(x,y), \\
&&G_i(\b(x),\b(y))=\b\circ G_i(x,y),\\
&&G_i(G_j(x,y), \a(z))=G_i(G_j(x,z),\a(y)),\\
&&G_i(\a\b(x),G_j(\a(y),z))-G_i(G_j(\b(x)\a(y)),\b(z))\nonumber\\
&& -G_i(\a\b(y),G_j(\a(x), z))+G_i(G_j(\b(y),\a(x)), \b(z))=0.
\end{eqnarray}

For $n=0$, this means  $A=A_0$ is a BiHom-Novikov algebra. For $n=1$, we obtain
some results for $G_1$:
\begin{eqnarray}
&&G_i(\a(x),\a(y))=\a\circ G_i(x,y), \\
&&G_i(\b(x),\b(y))=\b\circ G_i(x,y),\\
&&G_1(x\ast y, \a(z))+G_1(x, y)\ast \a(z)=G_1(x\ast z,\a(y))+G_1(x, z)\ast \a(y),\\
&&\a\b(x)\ast  G_1( \a(y), z)- G_1(\a(x),\a(y))\ast\b(z)-\a\b(y) \ast G_1( \a(x),  z)+G_1(\a(y)\ast\a(x))\ast \b(z) \nonumber\\
&&+\a\b(x)\ast G_1(\a(y),z)-G_1(\b(x), \a(y))\ast\b(z)-\a\b(y)\ast G_1(\a(x), z)+G_1(\b(y),\a(x))\ast \b(z)=0.\nonumber\\
&&\hspace{7cm}
\end{eqnarray}

 \noindent{\bf Example  4.5.}  By equation (2.8) with a fixed derivation
$D$, we obtain a  BiHom-Novikov algebra $(A, \star, \a, \b)$. Now we define a new multiplication $\ast_{\xi}$ on $A$ as
\begin{eqnarray}
x\ast_{\xi} y=\a(x)D(\b(y))+\xi xy.
\end{eqnarray}
for all $x,y\in A$, $\xi\in \mathbb{F}$. Then $(A, \star, \a, \b)$ is an infinitesimal deformation  BiHom-Novikov algebra,  we let $G_1(x,y)=\a(x)\b(y)$, and it is straightforward to check that  equations (4.10)-(4.13) hold.

Hence, we have the following result.

 \noindent{\bf Theorem  4.6.}  Let $(A, \mu)$ be a commutative algebra, $\a,\b: A\rightarrow A$ be two commuting algebra morphism,  and $D : A \rightarrow A$ be an even derivation such that $D\a=\a D, D\b=\b D$. Then $(A, \ast_{\xi}, \a, \b)$ is a BiHom-Novikov algebra,
where $\ast_{\xi}$ is defined as
\begin{eqnarray}
x\ast_{\xi} y=\a(x)D(\b(y))+\xi xy,
\end{eqnarray}
for all $x,y\in A$, $\xi\in \mathbb{K}$.

For two $\mathbb{K}$-bilinear maps $f,g:A\times A\rightarrow A$, define a map: $f\circ_{\a,\b}g: A[[t]]\times A[[t]]\times A[[t]]\rightarrow A[[t]] $ by
\begin{eqnarray*}
f\circ_{\a,\b}g(x,y,z)&=&f(\a\b(x),g(\a(y),z))-f(g(\b(x)\a(z)),\b(y))\\
&& -f(\a\b(y),g(\a(x), z))+f(g(\b(y),\a(x)), \b(z))¡£
\end{eqnarray*}
Using equations (4.8) and (4.9), we have
\begin{eqnarray*}
\sum_{i+j=n} G_i\circ_{\a,\b}G_j=0.
\end{eqnarray*}
For $n=1$, we have
\begin{eqnarray}
G_0\circ_{\a,\b}G_1+G_1\circ_{\a,\b}G_0=0.
\end{eqnarray}
For any $n\geq2$, we have
\begin{eqnarray}
-(G_0\circ_{\a,\b}G_n+G_n\circ_{\a,\b}G_0)=G_1\circ_{\a,\b}G_{n-1}+G_2\circ_{\a,\b}G_{n-2}+\cdot\cdot\cdot+G_{n-1}\circ_{\a,\b}G_0.~~~~
\end{eqnarray}

By (4.5)and (4.6) it follows that $G_i\in C_{\a,\b}^2(A,A)$. It can also be verified that $ G_i\circ_{\a,\b}G_j\in C_{\a,\b}^3(A,A)$.
In general, if $f,g\in C_{\a,\b}^2(A,A)$, then $ f\circ_{\a,\b}g\in C_{\a,\b}^3(A,A)$.  Note that the definition of
coboundary operator, which implies $\delta_{hom}^2G_n=G_0\circ_{\a,\b}G_n+G_n\circ_{\a,\b}G_0$ for $n=0,1,2\cdot\cdot\cdot$.
Hence (4.16) and (4.17) can be rewritten as
\begin{eqnarray*}
\delta_{hom}^2G_1=0,\\
-\delta_{hom}^2G_n=G_1\circ_{\a,\b}G_{n-1}+G_2\circ_{\a,\b}G_{n-2}+\cdot\cdot\cdot+G_{n-1}\circ_{\a,\b}G_0.
\end{eqnarray*}
Then $G_1$ is a 2-BiHom-cocycle.
\hfill $\square$

 \noindent{\bf Definition  4.7.} Let $(A,\ast, \a,\b)$ be a  regular BiHom-Novikov algebra. Suppose that $g_t(x, y)=G_0(\a^{-1}\b(x),y)+\sum_{i>0}G_i(x,y)t^i$
 and $g'_t(x, y)=\sum_{i\geq0}G'_i(x,y)t^i$ are two 1-parameter formal deformations of
$(A,\ast, \a,\b)$. They are called equivalent, denoted by $g_t\sim g'_t$, if there is a formal isomorphism
of $\mathbb{K}[[t]]$-modules
\begin{eqnarray*}
\phi_t(x)=\sum_{i\geq0}\phi_i(x)t^i: (A[[t]], g_t, \a,\b)\rightarrow (A[[t]], g'_t, \a,\b),
\end{eqnarray*}
where each  $\phi_i: T\rightarrow T$ is a $\mathbb{K}$-linear map and $\phi_0=id_A$, satisfying
\begin{eqnarray*}
\phi_t\circ \a=\a\circ \phi_t,\\
\phi_t\circ \b=\b\circ \phi_t,\\
\phi_t\circ g_t(x,y)=g'_t(\phi_t(x),\phi_t(y)).
\end{eqnarray*}
When $G_1=G_2=\cdot\cdot\cdot=0,g_t=G_0$ is said to be the null deformation. A 1-parameter for-
mal deformation $g_t$ is called trivial if $g_t\sim G_0$. A BiHom-Novikov algebra $(T,[\cdot,\cdot,\cdot], \a,\b)$
is called analytically rigid, if every 1-parameter formal deformation $g_t$ is trivial.

\noindent{\bf Theorem  4.8.}  Let  $g_t(x, y)=G_0(\a^{-1}\b(x),y)+\sum_{i>0}G_i(x,y)t^i$
 and $g'_t(x, y)=\sum_{i\geq0}G'_i(x,y)t^i$ be equivalent  1-parameter formal deformations of
$(A,\ast, \a,\b)$. Then $G_1$ and $G'_1$ belong to the
same cohomology class in $H_{\a,\b}^2(A,A)$.

{\bf Proof.} Suppose that $\phi_t(x)=\sum_{i\geq 0}\phi_i(x)t^i$ is the formal $\mathbb{K}[[t]]$-module  isomorphism
such that $\phi_t \circ \a=\a\circ \phi_t, \phi_t \circ \b=\b\circ \phi_t$ and
\begin{eqnarray*}
\sum_{i\geq 0}\phi_i(G_0(\a^{-1}\b(x),y)+\sum_{j>0}G_j(x,y)t^j)t^i=\sum_{i\geq 0}G'_{i}(\sum_{k\geq 0}\phi_k(x)t^k, \sum_{k\geq 0}\phi_l(y)t^l)t^i.
\end{eqnarray*}
It follows that
\begin{eqnarray*}
\sum_{i+j=n}\phi_i(G_0(\a^{-1}\b(x),y)+\sum_{j>0}G_j(x,y))t^{i+j}=\sum_{i+k+l=n}G'_{i}(\phi_k(x),\phi_l(y))t^{i+k+l}.
\end{eqnarray*}
In particular,
\begin{eqnarray*}
\sum_{i+j=1}\phi_i(G_0(\a^{-1}\b(x),y)+G_j(x,y)))t^{i+j}=\sum_{i+k+l=1}G'_{i}(\phi_k(x),\phi_l(y))t^{i+k+l}.
\end{eqnarray*}
that is,
\begin{eqnarray*}
G_1(x,y)+\phi_1(\a^{-1}\b(x)\ast y)=\phi_1(x)\ast y+x\ast \phi_1(y)+G'_1(x,y).
\end{eqnarray*}
Then $G_1-G'_1\in B_{\a,\b}^2(A,A)$.
\hfill $\square$

\noindent{\bf Theorem  4.9.} Suppose that $(A,\ast, \a,\b)$ is a regular BiHom-Novikov algebra such that $H_{\a,\b}^2(A,A)=0$. Then $(A,\ast, \a,\b)$ is analytically rigid.

{\bf Proof.} Let $g_t$ be a 1-parameter formal deformation of $(A,\ast, \a,\b)$. Suppose that $g_t=G_0+\sum_{i\geq n}G_it^i$. Then
\begin{eqnarray*}
\delta_{hom}^2G_n=G_1\circ_{\a,\b}G_{n-1}+G_2\circ_{\a,\b}G_{n-2}+\cdot\cdot\cdot+G_{n-1}\circ_{\a,\b}G_0=0,
\end{eqnarray*}
that is $G_n\in Z_{\a,\b}^2(A,A)=B_{\a,\b}^2(A,A)$. It follows that there exists $f_n\in C_{\a,\b}^1(A,A)$ such that $G_n=\delta_{hom}^1f_n$.

Let $\phi_t=id_A-f_nt^n: (A[[t]], g'_t, \a,\b)\rightarrow (A[[t]], g_t, \a,\b) $. Note that
\begin{eqnarray*}
\phi_t\circ \sum_{i\geq 0}f_n^it^{in}=\sum_{i\geq 0}f_n^it^{in}\circ \phi_t=id_{A[[t]]}.
\end{eqnarray*}
Then $\phi_t$ is a linear isomorphism. Moreover, $\phi_t\circ \a=\a\circ\phi_t $ and $\phi_t\circ \b=\b\circ\phi_t $.

Now we consider $g'_t(x,y)=\phi_t^{-1}g_t(\phi_t(x), \phi_t(y))$. It is straightforward to prove that $g'_t$ is a 1-parameter formal deformation of $(A,\ast, \a,\b)$ and $g_t\sim g'_t$, Suppose that $g'_t(x,y)=G'_0(\a^{-1}\b(x), y)+\sum_{i> 0}G'_i(x,y)t^i$. Then
\begin{eqnarray*}
&&(id_T-f_nt^n)(G'_0(\a^{-1}\b(x), y)+\sum_{i> 0}G'_i(x,y)t^i)\\
&=&  (G'_0(\a^{-1}\b(x), y)+\sum_{i\geq n}G'_i(x,y)t^i)(x-f_n(x)t^n)(y-f_n(y)t^n),
\end{eqnarray*}
i.e.,
\begin{eqnarray*}
&&G'_0(\a^{-1}\b(x), y)+\sum_{i> 0}G'_i(x,y)t^i-f_n(G'_0(\a^{-1}\b(x), y)+\sum_{i> 0}G'_i(x,y)t^i)\\
&=&\a^{-1}\b(x)\ast y-(f_n(x)\ast y+x\ast f_n(y))t^n+f_n(x)\ast f_n(y)t^{2n}\\
&&+\sum_{i\geq n}G_{i}(x,y)t^i-\sum_{i\geq n}(G_{i}(f_n(x),y)+G_{i}(x,f_n(y) )t^{i+n}+\sum_{i\geq n}G_i(f_n(x), f_n(y))t^{i+2n}.
\end{eqnarray*}
Then we have $G_1'=G_2'=...=G_{n-1}'=0$ and
\begin{eqnarray*}
G_n'(x,y)-f_n(\a^{-1}\b(x)\ast y) =-(f_n(x)\ast y+x\ast f_n(y))+G_n(x,y).
\end{eqnarray*}
Hence $G_n'=G_n-\delta^1f_n=0$ and $g'_t(x,y)=G'_0(\a^{-1}\b(x), y)+\sum_{i> 0}G'_i(x,y)t^i$. By induction, this procedure
ends with $g_t \sim G_0$, i.e., $(A,\ast, \a,\b)$. is analytically rigid. \hfill $\square$
\begin{center}
 {\bf ACKNOWLEDGEMENT}
 \end{center}

  The paper is partially supported by  the NSF of China (No. 11761017), the Fund of Science and Technology Department of Guizhou Province (No. [2016]1021),
 the Key University Science Research Project of Anhui Province (Nos. KJ2015A294, KJ2014A183 and KJ2016A545),
  the  NSF of Chuzhou University (Nos. 2014PY08 and 2015qd01).

\renewcommand{\refname}{REFERENCES}

\end{document}